\documentclass[12pt]{amsart}
\usepackage{graphicx}
\usepackage{amsmath}
\usepackage{amsfonts}
\usepackage{amssymb}
\usepackage{amscd}

\newcommand{\be}{\beta}

\newcommand{\Ga}{\Gamma}

\newcommand{\e}{\varepsilon}

\newcommand{\si}{\sigma}
\newcommand{\Si}{\Sigma}

\newcommand{\BQ}{\mathbb{Q}}
\newcommand{\BN}{\mathbb{N}}
\newcommand{\A}{\mathcal{A}}
\newcommand{\U}{\mathcal{U}}
\newcommand{\m}{\mathfrak{m}}
\newcommand{\ov}{\overline}

\newcommand{\val}{\mbox{val}_\beta}

\newcommand{\R}{\mathfrak{R}}

\newcommand{\wt}{\widetilde}

\newcommand{\0}{\mathbf 0}
\newcommand{\1}{\mathbf 1}

\renewcommand{\a}{\mathbf a}

\newcommand{\n}{\mathfrak{n}}

\renewcommand{\A}{\mathcal{A}}
\renewcommand{\phi}{\varphi}

\numberwithin{equation}{section}

\newtheorem{lemma}{Lemma}[section]

\newtheorem{prop}[lemma]{Proposition}
\newtheorem{thm}[lemma]{Theorem}

\theoremstyle{definition}
\newtheorem{Def}[lemma]{Definition}

\theoremstyle{remark}
\newtheorem{rmk}[lemma]{Remark}

\newtheorem{exam}[lemma]{Example}

\begin{document}

\title[Universal $\beta$-expansions]
{Universal $\beta$-expansions}
\author{Nikita Sidorov}
\address{Department of Mathematics, UMIST, P.O. Box 88,
Manchester M60 1QD, United Kingdom. E-mail:
Nikita.A.Sidorov@umist.ac.uk}
\date{\today}

\begin{abstract} Given $\be\in(1,2)$, a $\be$-expansion of a real
$x$ is a power series in base $\be$ with coefficients 0 and 1
whose sum equals $x$. The aim of this note is to study certain
problems related to the universality and combinatorics of
$\beta$-expansions. Our main result is that for any $\be\in(1,2)$
and a.e. $x\in (0,1)$ there always exists a universal
$\be$-expansion of $x$ in the sense of Erd\"os and Komornik, i.e.,
a $\be$-expansion whose complexity function is $2^n$. We also
study some questions related to the points having less than a full
branching continuum of $\be$-expansions and also normal
$\be$-expansions.
\end{abstract}

\maketitle

\section{Formulation of main results}\label{intro}

Let $\be\in(1,2)$ be our parameter, $\Si=\prod_1^\infty\{0,1\}$
and $x\ge0$. We will call a sequence $\e\in\Si$ a {\em
$\be$-expansion of $x$}, if it satisfies
\begin{equation}
x=\pi_\be(\e):=\sum_{n=1}^\infty\e_n\be^{-n}. \label{bexp}
\end{equation}
\begin{rmk}
Note that traditionally this term implies the greedy
$\be$-expansion of $x$ (see, e.g., \cite{Pa}) but for our purposes
it is better to use it in the above sense, because we will be
interested in {\bf all} $\be$-expansions of a given $x$. We hope
this will not cause any confusion.
\end{rmk}

It is clear that since $\e=(\e_n)_1^\infty\in\Si$, any $x$
representable in the form of the series~(\ref{bexp}), must belong
to the interval $I_\be:=[0,1/(\be-1)]$. On the other hand, each
$x\in I_\be$ does have at least one $\be$-expansion, namely, the
{\em greedy} $\be$-expansion: if $x\in[0,1)$, let
\[
T_\be(x)=\be x\bmod1,
\]
and put
\[
\e_n:=\lfloor\be T_\be^{n-1}x\rfloor,\quad n\ge1
\]
(here the power stands for the corresponding iteration,
$\lfloor\cdot\rfloor$ denotes the integral part of a number and
$\{\cdot\}$ stands for its fractional part). If $x\in[1,(\be-1))$,
then we put
\[
\ell=\min\,\{k\ge1 : x-\be^{-1}-\dots-\be^{-k}\in(0,1)\}
\]
and apply the greedy algorithm to $x-\be^{-1}-\dots-\be^{\ell}$ to
obtain the digits $\e_{\ell+1},\e_{\ell+2}$, etc. Finally, if
$x=1/(\be-1)$, then inevitably $\e_n\equiv1$.

An important property of the greedy $\be$-expansions consists in
their {\em monotonicity}, i.e., if $x<y$, then the greedy
$\be$-expansion of $x$ is lexicographically less than the one of
$y$. A detailed description of all possible greedy
$\be$-expansions for a given $\be$ was given by Parry \cite{Pa}
and is briefly described in Section~\ref{proofs}.

One of the intriguing questions regarding the $\be$-expansions
is as follows: does a given $x$ have $\be$-expansions different
from the greedy one, and if so, ``how many''? (cardinality,
dimension)

Recently the author proved the following metric result:

\begin{thm}\cite{Sid, AD}
For any $\be\in(1,2)$ a.e. $x\in I_\be$ has a continuum of
distinct $\be$-expansions.\label{cont}
\end{thm}

The proof given in \cite{Sid} is deliberately elementary;
however, in the survey paper \cite{AD} a more revealing
(dynamical) proof of this result is given. In the present paper
Theorem~\ref{cont} comes in a slightly stronger form -- see
Theorem~\ref{lessthancont}.

The main goal of this paper is to obtain a similar metric result
about universal $\be$-expansions which were introduced by Erd\"os
and Komornik \cite{EK98}. Recall their definition:

\begin{Def} A $\be$-expansion $\e=(\e_1,\e_2,\dots)$ is called
{\em universal} if for any finite 0-1 word $w$ there exists
$k\ge1$ such that $\e_k\dots\e_{k+N-1}=w$, where $N$ is the length
of $w$ (notation: $N=|w|$). In other words, its {\em complexity
function} must be $2^n$ (see, e.g., the survey paper \cite{Fer}
for the definition of complexity and dynamics-related results).
\end{Def}

In \cite{EK98} the authors concentrate their efforts mostly on the
case $x=1$; however, they present some results for general $x$ as
well. They show in particular that there exists $\be_0>1$ such
that for each $\be\in(1,\be_0)$ {\bf every} $x\in (0,1/(\be-1))$
has a universal $\be$-expansion. At the same time there exist
larger $\be$'s for which this is not the case (in particular,
$\be=\frac12(1+\sqrt5)$ -- see Counterexample below). The question
about the maximal possible $\be_0$ with this property remains
open; one of the obstacles proves to be the fact that if $\be$ is
a {\em Pisot number} (an algebraic integer greater than 1 whose
conjugates are all less than 1 in modulus), then, as was shown in
\cite{EK98}, $x=1$ cannot have a universal $\be$-expansion. It
would be conceivable if for {\em every}
$\be\in\left(1,\frac12(1+\sqrt5)\right)$ every $x$ had a universal
$\be$-expansion unless $\be$ is a Pisot number.

The main result of the present paper is

\begin{thm} For every fixed $\be\in(1,2)$ a.e. $x\in I_\be$ has a
universal $\be$-expansion.\label{1}
\end{thm}

It is natural to ask the question whether Theorems~\ref{cont} and
\ref{1} are related. The answer is negative in one direction:
namely, there exist $\be$ and $x$ having a continuum of distinct
$\be$-expansions, none of which is universal:

\smallskip\noindent\textbf{Counterexample.} Let
$\be=\frac12(1+\sqrt5)$ and $x=1/2\be$. As was shown by Vershik
and the author in \cite{SV}, the space of all $\be$-expansions
for $x$ is $\prod_1^\infty\{100,011\}$. It is easy to see that
the word 1010 cannot occur in any of them, i.e., none of them is
universal.

\smallskip
We believe (but have failed to show) that if $x$ has a universal
$\be$-expansion, then it has a continuum of distinct
$\be$-expansions. Since it is easy to show that any $x$ that has a
finite number of $\be$-expansions can be excluded (see
Remark~\ref{fini} below), the problem may be reformulated as
follows:

\smallskip\noindent\textbf{Open problem.}\textit{
Prove or disprove that there exist $\be\in(1,2)$ and $x$ which has
precisely $\aleph_0$ different $\be$-expansions such that one of
them is universal}.

Section~\ref{comb} contains a claim that improves
Theorem~\ref{cont}; more precisely, we deal with the set of points
$x$ whose {\em branching compactum} is not full (for example, such
are points which have not more than countable set of
$\be$-expansions). We show that this set is in a way very close to
the set of unique $\be$-expansions studied in \cite{GS} (see
Proposition~\ref{str} below) and in particular, has the Hausdorff
dimension strictly less than~1.

Finally, in Section~\ref{normal} we discuss more delicate
questions related to {\em normal} $\be$-expansions (see the
definition at the beginning of Section~\ref{normal}).

\section{Universality}\label{proofs}

This section is devoted completely to the proof of
Theorem~\ref{1}. Our method may be called {\em
anti-normalization}; let us recall first some definitions and
present a model example, for which our proof will be especially
simple and revealing.

\subsection{Necessary definitions}\label{defs} Let the sequence
$(a_n)_1^\infty$ be defined as follows: let
$1=\sum_{1}^{\infty}a_k' \be^{-k}$ be the greedy expansion of 1,
i.e, $a_n'=\lfloor\be T_\be^{n-1}1\rfloor,\ n\ge1$. If the tail
of the sequence $(a_n')$ differs from $0^\infty$, then we put
$a_n\equiv a_n'$. Otherwise let $k=\max\,\{j:a_j'>0\}$, and
$(a_1,a_2,\dots):= (a_1',\dots,a_{k-1}',a_k'-1)^\infty$. In the
seminal paper \cite{Pa} it is shown that for each greedy
expansion $\e$ in base $\be,\ (\e_n,\e_{n+1},\dots)$ is
lexicographically less (notation: $\prec$) than
$(a_1,a_2,\dots)$ for every $n\in\BN$. Moreover, it was shown
that, conversely, every sequence with this property is actually
the greedy expansion in base $\be$ for some $x\in [0,1)$.

Thus, the set of all greedy $\be$-expansions for $x\in[0,1)$ is
the {\em $\be$-compactum}
\begin{equation}
X_\be=\left\{\e\in\Si\mid
(\e_n,\e_{n+1},\dots)\prec(a_1,a_2,\dots),\
n\in\BN\right\}\label{Xbe}
\end{equation}
The sequences from the $\be$-compactum will be called {\em
$\be$-admissible} (or simply {\em admissible} if it is clear which
$\be$ is in question).

As is shown in \cite{Pa}, the map $\phi_\be:X_\be\to[0,1)$ defined
by the formula
\begin{equation}
\phi_\be(\e)=\sum_{n=1}^\infty\e_n\be^{-n}\label{bexp2},
\end{equation}
is in fact one-to-one with the exception of a countable set of
sequences.

Let $\tau_\be$ denote the shift on $X_\be$ (it is obvious from
(\ref{Xbe}) that the $\be$-compactum is shift-invariant). Then
$T_\be=\phi_\be\tau_\be\phi_\be^{-1}$, and it is shown in
\cite{Re} that there exists unique ergodic shift-invariant
probability measure $\nu_\be$ on $X_\be$ such that
$\phi_\be(\nu_\be)$ is equivalent to the Lebesgue measure. This
measure is positive on all cylinders $[\e_1=i_1,\dots,\e_n=i_n]$
provided $(i_1,\dots,i_n,0^\infty)$ is $\be$-admissible.

Define the {\em $\be$-value} of a 0-1 word $w=(x_k)_1^N$
($N\le\infty$) as follows:
$$
\val(w):=\sum_1^N x_k\be^{-k}.
$$
Let $\Si_\be$ denote the set of 0-1 sequences whose $\be$-value is
less than or equal to 1. The {\em normalization} (in base $\be$)
is the map from $\Si_\be$ to $X_\be$ defined by the formula
\begin{equation}
\n_\be:=\phi_\be^{-1}\circ\pi_{\be}, \label{norm}
\end{equation}
where $\pi_\be$ is given by (\ref{bexp}) and $\phi_\be$ is given
by (\ref{bexp2}). We will identify sequences whose tail is
$0^\infty$ with the corresponding finite words. So, the
normalization of a finite word can be finite (see below). For
more detail and ``automatic'' properties of normalization see
\cite{Fr}.

Thus, the operation of normalization assigns to each 0-1
sequence the admissible sequence with the same $\be$-value.

\begin{Def} Two 0-1 words of the same length (finite or
infinite) will be called {\rm equivalent}, if they have the same
$\be$-value.
\end{Def}

\subsection{Special case of Theorem~\ref{1}} Consider the case
$\be=G:= \frac12(1+\sqrt5)$. Let $W$ denote the set of all finite
0-1 words whose $\be$-value is less than 1, and $V$ denote the set
of all $\be$-admissible 0-1 words (obviously, $V\subsetneq W$).
Note that in this case admissibility simply means that there are
no two consecutive unities. As is well known, in this case the
normalization of a finite $w\in W$ is a finite word of the same
length (see, e.g., \cite{SV}).

Since the shift $(X_\be,\nu_\be,\tau_\be)$ is ergodic, for any
fixed $v\in V$, $\nu_G$-a.e. sequence $\e\in X_G$ contains $v$
infinitely many times.\footnote{More precisely, there exist
infinitely many $k$ such that $\e_k\dots\e_{k+|v|-1}=v$.} Let the
equivalence class of $v$ be $[v]=\{w_1,\dots,w_\ell\}$.

We perform the anti-normalization for $v$ as follows: take a
generic sequence $\e$ in question; first time we hit $v$, we
change it to $w_1$, next time we hit it -- to $w_2$, etc., until
we get to $w_\ell$. Note that this operation has not changed the
$\be$-value of $\e$. Now we fix some ordering of $V$, say, the
lexicographic one:
$V=\{v_1,v_2,\dots\}=\{0,1,00,01,10,000,001,\dots\}$ and perform
this operation consecutively for $v_1,v_2$, etc. Since each word
from $V$ occurs infinitely many times, we can avoid ``overlaps".
The resulting ``anti-normalized" sequence $\e'$ a
$\be$-expansion of a Lebesgue-generic $x$, and by our
construction, it contains all 0-1 words, i.e., is a universal
sequence.

The actual reason why the case $\be=G$ is so easy to deal with, is
the fact that $G$ is a {\em finitary} Pisot number, i.e., the
normalization of each finite word in base $G$ is finite as well.
For an arbitrary (even Pisot) $\be$ this is, generally speaking,
not true, and we will need a more delicate argument.

\subsection{General case} Note first that it suffices to prove
Theorem~\ref{1} only for $x\in(0,1)$. This is because for
$x\in(1,1/(\be-1))$ one can find $\ell\ge1$ such that
$y=x-\be^{-1}-\dots-\be^{-\ell}\in[0,1)$, then apply the theorem
to $y$ (i.e., a generic $y$ has a universal $\be$-expansion
$(y_1,y_2,\dots)$). Finally, $(1,\dots,1,y_1,y_2,\dots)$ (with
$\ell$ unities) is a universal $\be$-expansion of $x$.

Since the proof is somewhat technical, we would like to present
a sketch first and then fill up the details. Let $V, W$ be as in
the previous subsection, $w\in W$; then there exists $v\in V$
whose $\be$-value is slightly greater than the $\be$-value of
$w$ (just consider $v'=\n_\be(w)$, replace 0 by 1 at a
sufficiently large coordinate of $v'$ and drop the rest of it).

A generic sequence $\e\in X_\be$ is of the form
$(\dots,v,(tail))$ and we ``anti-normalize" it into
$(\dots,w,(tail'))$, where the dots denote one and the same
symbols, $tail'=(\e_n',\e_{n+1}',\dots)$ is admissible, and
$tail'=v+tail-w$ (we identify an admissible sequence with its
$\be$-value). Since we have chosen $v$ slightly greater than
$w$, the $\be$-value of $(\e_n',\e_{n+1}',\dots)$ fills some
interval $(a,b)\subset (0,1)$, where $a=a(w), b=b(w)$. Since
$tail$ is ``random", so is $tail'$ (more precisely, its shift
$(\e_n',\e_{n+1}',\dots)$). Hence we can repeat this procedure
{\em ad infinum}; the claim follows from the fact that for a.e.
$x\in (a,b)$, its greedy $\be$-expansion contains every
admissible word infinitely many times (this is a trivial
consequence of Poincar\'e's recurrence theorem for the shift
$\tau_\be$).

To turn this sketch into a real proof, we have to clarify the
following points:
\begin{enumerate}
\item accurate choice of $v$;
\item ``randomness" of $tail'$.
\end{enumerate}

\medskip\noindent (1) Let $w=x_1\dots x_k$ and
$v'=(\e_1,\dots,\e_k,\dots)=\n_\be(w)$. Put
$$
n:=\min\,\{j\ge1 : \e_{k+j}=0\,\,\mbox{and}\,\,
\e_1\dots\e_{k+j-1}1\,\,\mbox{is\ $\be$-admissible}\}.
$$
The number $n$ is well defined, because the tail of $v'$ does
not coincide with the tail of $(a_i)_1^\infty$ (see the
beginning of the section), whence one can always increase $v'$
at a sufficiently large coordinate.

Put $v:=\e_1\dots\e_{k+n-1}1$. Now $a:=\val(v)-\val(w)\in(0,1)$
(it is positive by the monotonicity of the greedy
$\be$-expansions -- see Section~\ref{intro}). To determine $b$,
we consider the sequence $\wt v$ which is defined as the largest
possible $\be$-admissible sequence beginning with $v$. Then
$b:=\val(\wt v)-\val(w)<1$.

\smallskip\noindent (2) Put
\begin{align*}
E_w^{(j)}:=\{&x\in(0,1) : \e=(\e_1,\dots)\,\mbox{is the greedy
$\be$-expansion of $x$ and}\\ &\e_{j+1}\dots\e_{j+|v|}=v\},
\end{align*}
(here $v=v(w)$ as above) and $E_w:=\cup_{j\ge0}E_w^{(j)}$.
Obviously, $\mathcal L(E_w)=1$, where $\mathcal L$ denotes
Lebesgue measure. The relation $v+tail=w+tail'$ can be rewritten
in the following way: we have for $x\in E_w^{(j)}$:
\[
\val(v)+\be^{-|v|}T_\be^{j+|v|}x=\val(w)+\be^{-|w|}y,
\]
whence
\[
y=c_1(w)+c_2(w)T_\be^{j+|v|}x\in (a(w),b(w)),
\]
where $c_1(w), c_2(w)$ are some constants. Hence in view of
$T_\be(\mathcal L)$ being equivalent to $\mathcal L$
(\cite{Re}), Lebesgue measure of all possible $y$'s in
$(a(w),b(w))$ is full. Therefore, a generic $x$ leads to a
generic $y$, and we may repeat this operation for all $w$, thus
constructing a universal $\be$-expansion of a generic
$x$.\footnote{Note that from the proof it follows that $y$ being
in $(a,b)$ (whereas $x$ could assume any value in $(0,1)$) does
not affect the choice of the next interval $(a,b)$.}

Theorem~\ref{1} is proved.

\section{Combinatorics and branching}\label{comb}

\subsection{Unique expansions}
We need first to recall some facts about unique $\be$-expansions.
Namely, $x\in I_\be$ will be said to have {\em unique}
$\be$-expansion if the greedy $\be$-expansion is the only one
$\be$-expansion for $x$. Let $\A_\be$ denote the set of such
$x\in(0,1/(\be-1))$'s ($x=0$ and $x=1/(\be-1)$ obviously have a
unique $\be$-expansion). It is shown in \cite{EJK} that if
$\be<G$, then $\A_\be=\emptyset$. A natural question to ask is
about its properties when $G\le\be<2$. The following theorem has
been recently proved by P.~Glendinning and the author:

\begin{thm} \label{uniq} \cite{GS} The set $\A_\be$ has measure
zero for any $\be\in(1,2)$. The cardinality of the set $\A_\be$
is
\begin{itemize}
\item [(i)] $\aleph_0$ if $\be\in(G,\be_*)$ and
\item [(ii)] $2^{\aleph_0}$ if $\be\in[\be_*,2)$.
\end{itemize}
Moreover, if $\be=\be_*$, then $\A_\be$ is a Cantor set of zero
Hausdorff dimension, and if $\be\in(\be_*,2)$, then
$0<\dim_H(\A_\be)<1$.
\end{thm}

Here $\be_*=1.787231650\dots$ is the {\em Komornik-Loreti
constant}, i.e., the smallest $\be$ such that $x=1$ has a unique
$\be$-expansion. In \cite{KL} it is shown that in fact $\be_*$
is the unique solution of the equation
$$
\sum_{n=1}^\infty \m_n x^{-n+1}=1,
$$
where $\m=(\m_n)_1^\infty$ is the Thue-Morse sequence
\cite{AllSh}:
\begin{equation*}
\m=0110\,\,1001\,\,1001\,\,0110\,\,1001\,\,0110\,\,0110\,\,1001
\dots
\end{equation*}
In \cite{GS} we have given a symbolic description of unique
$\be$-expansions. Namely, let $\mathbf a=(a_1,a_2,\dots)$ (see
Section~\ref{proofs}), and $\si$ denote the shift on $\Si$.
Define
\[
\U_\be:=\{\e\in\Si: \overline\a\prec \si^n\e\prec\a, \ n\ge0 \},
\]
where bar denotes the {\em inversion}, i.e., $\ov0=1,\ov1=0$.

In \cite{GS} it is shown that any unique $\be$-expansion which is
neither $0^\infty$ nor $1^\infty$, is of the form $0^s\e$ or
$1^s\e$, where $s\ge0$, and $\e\in\U_\be$.

\subsection{``Less than continuum" of $\be$-expansions} We will
show that having ``less than the full continuum" of possible
$\be$-expansions is almost the same cardinality-wise as having a
unique one. Let $\R_\be(x)$ denote the set of all
$\be$-expansions of $x$. We are going to construct by induction
the {\em branching compactum}
$\Ga_\be(x)\subset\prod_1^\infty\{\0,\1\}$.

Firstly, if $x\in\A_\be$, then we define
$\Ga_\be(x):=\{\0^\infty\}$; otherwise, there exists a {\em
branching}, i.e., there exist
$(\e_1,\dots,\e_n,\e_{n+1},\dots)\in\R_\be(x)$ and
$(\e_1,\dots,\e_n,\e'_{n+1},\e'_{n+2}\dots)\in\R_\be(x)$ for
some $n\ge0$, and $\e_{n+1}\neq\e_{n+1}'$. Thus, we can make a
choice for the first symbol in $\Ga_\be(x)$: it is $\1$ if we
choose the ``lower branch" (i.e., zero at the $n$'th place) and
$\0$ otherwise.\footnote{This in fact corresponds to the
dynamical model described in detail in \cite[\S2]{AD}.}
Performing the same operation for $(\e_{n+1},\e_{n+2},\dots)$
and $(\e'_{n+1},\e'_{n+2},\dots)$ yields the second symbol in
$\Ga_\be(x)$, etc. If one of the ``tails" happens to be a unique
$\be$-expansion, we assume for simplicity that the remaining
symbols are all $\0$'s.

\begin{exam} For $\be=G$ and $x=\be^{-1}$, as is well known,
$\R_\be(x)=\{10^\infty,0110^\infty,010110^\infty,\dots\}$,
whence
$\Ga_\be(x)=\{\0^\infty,\linebreak[0]\1\0^\infty,\linebreak[0]
\0\1\0^\infty,\linebreak\0\0\1\0^\infty\dots\}$. On the other
hand, in Counterexample described in Section~\ref{intro},
$\Ga_\be(x)=\prod_1^\infty\{\0,\1\}$.
\end{exam}

By our construction, to every $\be$-expansion of $x$ one can
assign the (unique) sequence from $\Ga_\be(x)$, i.e.,
$\R_\be(x)$ and $\Ga_\be(x)$ are naturally isomorphic. Put
$$
\mathcal C_\be:=\Bigl\{x\in (0,1/(\be-1)) :
\Ga_\be(x)\neq\prod_1^\infty\{\0,\1\}\Bigr\}.
$$
Thus, the set of $x$'s, for which
$\mbox{card}\,\R_\be(x)<2^{\aleph_0}$, is a subset of $\mathcal
C_\be$. Note that in \cite[Theorem~2.18]{AD} we have in fact
shown that for every $\be$ the set $\mathcal C_\be$ has zero
Lebesgue measure. Here we would like to make this result more
precise.

The following auxiliary claim is straightforward:

\begin{lemma}\label{shift}
$x\in\mathcal C_\be$ if and only if there exists its
$\be$-expansion $\e$ and $n\ge0$ such that $\si^n\e$ is a unique
$\be$-expansion.
\end{lemma}
This simple observation helps us to refine Theorem~\ref{cont}.

\begin{lemma}\label{count21}
There exists a map $\psi_\be:\mathcal C_\be\to\A_\be$ which is
countable-to-one.
\end{lemma}
\begin{proof}
By the above, if $x\in\mathcal C_\be$, then
$x=\sum_1^{n-1}\e_j\be^{-j}+\be^{-n}y$, where $\e_j\in\{0,1\}$ and
$y\in\A_\be$. We define the map $\psi_\be: x\mapsto y$. The choice
of a specific $y$ is unimportant; for example, if there multiple
$y$'s, choose the smallest $n$ first, and if there is still a
choice, choose the lexicographically smallest
$(\e_1,\dots,\e_{n-1})$.

Now, if we have also $x'=\sum_1^{n'-1}\e'_j\be^{-j}+\be^{-n'}y$,
then $x-\be^k x'\in\BQ(\be)$ for $k=n'-n$, whence for a given
$x\in\psi_\be^{-1}\{y\}$ there can be not more than a countable
set of $x'$'s from the same preimage.
\end{proof}

Recall now that in \cite[Theorem~3]{EJK} quoted above, it was in
fact shown that for any $\be\in(1,G),\
\mbox{card}\,\R_\be(x)=2^{\aleph_0}$ for {\em every} $x\in
(0,1/(\be-1))$.\footnote{Actually, from their proof it even
follows that $\mathcal C_\be=\emptyset$ for every $\be<G$.} This
result is in a way best possible, because for $\be=G$ there is
already a countable set of points, each of which has $\aleph_0$
$\be$-expansions (for instance, $x=1$), and for $\be>G$, as we
know, there are points which even have a unique $\be$-expansion.

Nevertheless, some improvement of Theorem~\ref{cont} for $\be>G$
is possible. Namely, we show that having a non-full branching is
very close to having just a single $\be$-expansion.

Let $w_n:=\m_2\dots\m_{2^n+1},\ n\ge0$, where $\m$ is the
Thue-Morse sequence. That is, $w_0=1,w_1=11,w_2=1101,w_3=1101\
0011$, etc.

\begin{prop}\label{str} {\em(1)} For any $\be\in(G,\be_*)$ we have
$\mathcal C_\be\subset\BQ(\be)$. More precisely, every
$x\in\mathcal C_\be$ has an eventually periodic $\be$-expansion
with the period $0^\infty$, $1^\infty$ or $w_n\ov w_n$ for some
$n\ge0$.\newline {\em (2)} For $\be\in[\be_*,2)$,
\[
\dim_H \mathcal C_\be=\dim_H \A_\be\in[0,1).
\]
\end{prop}
\begin{proof} (1) By \cite[Proposition~13]{GS}, $\U_\be$ contains
only eventually periodic sequences with the period $w_n\ov w_n$
for some $n\ge0$ if $\be\in(G,\be_*)$. Now the claim follows
directly from Lemma~\ref{shift}. \newline (2) is a consequence
of Lemma~\ref{count21}.
\end{proof}

As a corollary of \cite[Theorem~3]{EJK} and
Proposition~\ref{str} we obtain

\begin{thm}\label{lessthancont}
The set
\[
\left\{x\in(0,1/(\be-1)) :
\mathrm{card}\,\R_\be(x)<2^{\aleph_0}\right\}
\]
is
\begin{itemize}
\item empty if $\be\in (1,G)$;
\item a proper subset of $\BQ(\be)$ if $G\le\be<\be_*$;
\item a continuum of Hausdorff dimension 0 if $\be=\be_*$;
\item a continuum of Hausdorff dimension strictly between 0
and 1 if $\be\in(\be_*,2)$.
\end{itemize}
\end{thm}

\begin{rmk}\label{fini}
Note that if $\R_\be(x)$ is finite, then for every
$\e\in\R_\be(x)$ there exists $n=n(\e)\ge0$ such that $\si^n\e$
is unique. Hence such a sequence cannot be universal, because a
unique $\be$-expansion cannot contain, for instance, the word
$0^s$ for $s$ large enough.

On the other hand, if $\R_\be(x)$ is countable, there will be both
sequences whose $n$'th shift is unique but also inevitably those
not having this property. This is the main obstacle for an easy
solution of the open problem mentioned in the end of
Section~\ref{intro}.
\end{rmk}

\section{Normal $\be$-expansions}\label{normal}

We know from Theorem~\ref{1} that for a given $\be$ a.e. $x$ has
at least one universal $\be$-expansion. Note that by the
ergodicity of the shift $\tau_\be$, each admissible block occurs
in the greedy expansion of a generic $x$ with a positive limiting
frequency. Thus, the proof given in Section~\ref{proofs} can be
easily modified to yield

\begin{prop}
Given $\be\in(1,2)$, a.e. $x\in I_\be$ has a universal
$\be$-expansion with a positive limiting frequency of each 0-1
block.
\end{prop}

It is thus natural to ask the following question: is it true that
for every $\be$ a.e. $x\in I_\be$ has a \textit{normal}
$\be$-expansion, i.e., the one for which the limiting frequency of
each 0-1 block $B$ is exactly $2^{-|B|}$?\footnote{Yet again, we
hope there will be no confusion with the notion of normal {\em
greedy} $\be$-expansions -- see, e.g., \cite{BM}.} A partial
answer to this question is

\begin{thm}
There exists a set $E\subset(1,2)$ of full Lebesgue measure such
that for each $\be\in E$, Lebesgue-a.e. $x\in I_\be$ has a
normal $\be$-expansion.
\end{thm}
\begin{proof}
Let $p$ denote the product measure
$\prod_1^\infty\{\frac12,\frac12\}$ on $\Si$, and
$\mu_\be=\pi_\be(p)$, where $\pi_\be$ is given by (\ref{bexp}).
This measure is called the {\em Bernoulli convolution}
parameterized by $\be$ (see, e.g., \cite{PSS}).

Note first that the claim in question is valid for every $\be$
and $\mu_\be$-a.e. $x$ -- it suffices to consider a set
$\mathfrak N$ of ``normal" sequences in $\Si$ (which by the SLLN
has $p$-measure~1) and take $\mathfrak N_\be:=\pi_\be(\mathfrak
N)$. This set will have full $\mu_\be$-measure, and clearly,
every $x\in\mathfrak N_\be$ has a normal $\be$-expansion. To end
the proof of the theorem, it suffices to recall that by the
famous theorem due to B.~Solomyak \cite{So}, for a.e. $\be$ the
Bernoulli convolution $\mu_\be$ is absolutely continuous with
respect to the Lebesgue measure on $I_\be$, whence for a.e.
$\be$ the set $\mathfrak N_\be$ has Lebesgue measure~1 as well.
\end{proof}

\begin{rmk} We believe $\mathfrak N_\be$ has Lebesgue
measure~1 for all $\be$, even if $\be$ is a Pisot number (it is
well known that $\mu_\be$ in this case is singular \cite{E}). We
plan to return to this problem in the future.
\end{rmk}

\medskip\noindent \textbf{Acknowledgment.} The author's research
was supported by the \linebreak EPSRC grant no GR/R61451/01. The
author is grateful to Vilmos Komornik for stimulating
discussions.

\end{document}